\documentclass[12pt]{amsart}
\usepackage{amsmath,amsthm,enumerate,graphics,amsfonts,latexsym,amsopn,verbatim,amscd,amssymb}
\usepackage{color}
\usepackage{tikz}

\theoremstyle{plain}
\newtheorem{thm}{Theorem}[section]

\newtheorem{cor}[thm]{Corollary}

\theoremstyle{definition}
\newtheorem{rem}[thm]{Remark}

\DeclareMathOperator{\cnspec}{CN-spec}

\DeclareMathOperator{\CN}{CN}
\DeclareMathOperator{\Ecn}{E_{CN}}

\begin{document}

\title[CN-spectrum and energy of commuting conjugacy class graph]{Common neighbourhood spectrum and energy of commuting conjugacy class graph}

\author[F. E. Jannat and R. K. Nath ]{Firdous Ee Jannat and Rajat  Kanti Nath}

\address{Firdous Ee Jannat,  Department of Mathematical Science, Tezpur  University, Napaam -784028, Sonitpur, Assam, India.}

\email{firdusej@gmail.com}

\address{Rajat  Kanti Nath, Department of Mathematical Science, Tezpur  University, Napaam -784028, Sonitpur, Assam, India.}

\email{rajatkantinath@yahoo.com}

\begin{abstract} 
In this paper we compute  common neighbourhood (abbreviated as CN) spectrum and energy of commuting conjugacy class graph of several families of finite non-abelian groups.  
As a consequence of our results we show that the commuting conjugacy class graphs of the groups $D_{2n}$, $T_{4n}$, $SD_{8n}$, $U_{(n,m)}$, $U_{6n}$, $V_{8n}$,  $G(p, m, n)$ and  some  families of groups whose central quotient is isomorphic to $D_{2n}$ or $Z_p \times Z_p$, for some prime $p$, are CN-integral but not CN-hyperenergetic.
\end{abstract}


\maketitle

\section{Introduction} \label{S:intro} 

Let $G$ be a finite non-abelian group. The commuting conjugacy class graph of $G$, denoted by $\Gamma(G)$, is defined as a simple undirected graph whose vertex set is the set of conjugacy classes of the non-central elements of $G$ and two vertices $a^G$ and $b^G$ are adjacent if there exists some elements $a'\in a^G$ and $b'\in b^G$ such that $a'b'=b'a'$. The study of commuting conjugacy class graphs of groups was initiated by  Herzog, Longobardi and Maj \cite{HLM} in the year 2009. In 2016, Mohammadian et al.\cite{MEFW-2016} have characterized finite groups such that their commuting conjugacy class graph is triangle-free.  
Later on Salahshour and Ashrafi \cite{SA-2020,SA-CA-2020}, obtained structures of commuting conjugacy class graph of several families of finite CA-groups.  Salahshour \cite{Salah-2020} also described $\Gamma(G)$ for the groups whose central quotient is isomorphic to a dihedral group.
In \cite{BN-2021}, Bhowal and Nath have characterized certain finite groups such that   $\Gamma(G)$ is hyperenergetic, L-hyperenergetic and Q-hyperenergetic.
They also have characterized certain finite groups such that $\Gamma(G)$ is planar, toroidal double-toridal and triple-toroidal in  \cite{BN-2022}.

\par Let $\mathcal{G}$ be a simple graph with vertex set $V(\mathcal{G}) := \{v_i:i=1,2,\ldots,n\}$. The common neighbourhood (abbreviated as CN) of two distinct vertices $v_i$ and $v_j$, denoted by $C(v_i,v_j)$ is the set of all vertices other than $v_i$ and $v_j$, which are adjacent to both $v_i$ and $v_j$. The common neighbourhood matrix of $\mathcal{G}$, denoted by $\CN(\mathcal{G})$, is defined as 
$$(\CN(\mathcal{G}))_{i,j}=\begin{cases}
|C(v_i,v_j)|, &\text{ if } i\neq j\\
0,&\text{ if } i=j.
\end{cases}$$
 The set of all eigenvalues of $\CN(\mathcal{G})$ with multiplicities, denoted by $\cnspec(\mathcal{G})$, is called the common neighbourhood spectrum (abbreviated as CN-spectrum) of $\mathcal{G}$. If $\lambda_1,\lambda_2,\ldots,\lambda_k$ are the eigenvalues of $\CN(\mathcal{G})$ with multiplicities  $\alpha_1,\alpha_2,\ldots,\alpha_k$ respectively then we write $\cnspec(\mathcal{G})=\{\lambda_1^{\alpha_1},\lambda_2^{\alpha_2},\ldots,\lambda_k^{\alpha_k}\}$.  A graph $\mathcal{G}$ is called CN-integral if $\cnspec(\mathcal{G})$ contains only integers. The common neighbourhood energy (abbreviated as CN-energy) of a graph $\mathcal{G}$, denoted by $\Ecn(\mathcal{G})$, is defined as
$$\Ecn(\mathcal{G})=\sum_{i=1}^n\alpha_i|\lambda_i|.$$
 A graph $\mathcal{G}$ is called CN-hyperenergetic if $\Ecn(\mathcal{G})>\Ecn(K_n)$, where $n=|V(\mathcal{G})|$. If $\Ecn(\mathcal{G})=\Ecn(K_n)$ then $\mathcal{G}$ is called  CN-borderenergetic.  
  In 2011, Alwardi,  Soner and  Gutman \cite{ASG}  introduced this concepts of CN-spectrum and CN-energy of a graph. 
 
In this paper we compute  common neighbourhood  spectrum and energy of commuting conjugacy class graph of several families of finite non-abelian groups.  As a consequence of our results we show that the commuting conjugacy class graphs of the groups $D_{2n}$, $T_{4n}$, $SD_{8n}$, $U_{(n,m)}$, $U_{6n}$, $V_{8n}$,  $G(p, m, n)$ and  some  families of groups whose central quotient is isomorphic to $D_{2n}$ or $Z_p \times Z_p$, for some prime $p$, are CN-integral but not CN-hyperenergetic.



\section{Computation of CN-spectrum and CN-energy}

In this section we compute  CN-spectrum and CN-energy of commuting conjugacy class graph of several families of finite groups.
We write $\mathcal{G}=\mathcal{G}_1\cup \mathcal{G}_2\cup\ldots\cup\mathcal{G}_n$ to denote that a graph $\mathcal{G}$ has $n$ components namely $\mathcal{G}_1$, $\mathcal{G}_2$,\ldots,$\mathcal{G}_n$. Also $lK_m$ denotes the disjoint union of $l$ copies of the complete graph $K_m$ on $m$ vertices. By \cite[Theorem 1]{FSN-2021} and \cite[Theorem 2.3]{NFDS-2021} we get the following result which is very useful in computing CN-spectrum and CN-energy of commuting conjugacy class graph of the groups considered in this paper. 

%


\begin{thm}\label{CNS2}
	Let $\mathcal{G} = l_1K_{m_1} \cup l_2K_{m_2}\cup l_3K_{m_3}$, where $l_iK_{m_i}$ denotes the disjoint union of $l_i$ copies of the complete graphs $K_{m_i}$ on ${m_i}$ vertices for $i = 1, 2, 3$. Then
\begin{align*}
\cnspec(\mathcal{G})=\{(-&(m_1-2))^{l_1(m_1-1)}, ((m_1-1)(m_1-2))^{l_1}, (-(m_2-2))^{l_2(m_2-1)},\\
& ((m_2-1)(m_2-2))^{l_2}, (-(m_3-2))^{l_3(m_3-1)},	((m_3-1)(m_3-2))^{l_3}\}
\end{align*}
and
$$
\Ecn(\mathcal{G})=2l_1(m_1-1)(m_1-2) + 2l_2(m_2-1)(m_2-2) + 2l_3(m_3-1)(m_3-2).
$$
\end{thm}

\subsection{Certain $2$-generated finite groups}
In this subsection we consider the dihedral groups, dicyclic group, semidihedral group along with some other $2$-generated finite group and compute CN-spectrum and CN-energy of their  commuting conjugacy class graph. 
\begin{thm}\label{Thm 3.1}
If $G$ is the dihedral group $D_{2n}=\langle x,y:~x^{2n}=y^2=1,~yxy^{-1}=x^{-1} \rangle$ then CN-spectrum and CN-energy of $\Gamma(G)$ are given by 
$$\cnspec(\Gamma(G))=\begin{cases}
\{(-\frac{1}{2}(n-5))^{\frac{1}{2}(n-3)},(\frac{1}{4}(n-3)(n-5))^1,0^1\},~~~\text{ if }~ 2\nmid n \\
\{(-\frac{1}{2}(n-6))^{\frac{1}{2}(n-4)},(\frac{1}{4}(n-4)(n-6))^1,0^2\},~~~\text{ if }~  2\mid n
\end{cases}
$$
and
$$
\Ecn(\Gamma(G))=\begin{cases}
\frac{1}{2}(n-3)(n-5),~~~\text{ if }~ 2\nmid n\\
\frac{1}{2}(n-4)(n-6),~~~\text{ if }~  2\mid n.
\end{cases}
$$

\end{thm}
\begin{proof}
By \cite[Proposition 2.1]{SA-CA-2020}, we have
$$\Gamma(G)=\Gamma(D_{2n})=\begin{cases}
K_{\frac{n-1}{2}}\cup K_1, & \text{ if }~  2\nmid n \\
K_{\frac{n}{2}-1}\cup 2K_1, & \text{ if }~  2\mid n \text{ and } \frac{n}{2} \text{ is even } \\
 K_{\frac{n}{2}-1}\cup K_2, & \text{ if }~  2\mid n \text{ and } \frac{n}{2} \text{ is odd.}
\end{cases}$$
Now by applying Theorem \ref{CNS2}, we get
$$
\cnspec(\Gamma(G))=\begin{cases}
\{(-(\frac{n-1}{2}-2))^{(\frac{n-1}{2}-1)},((\frac{n-1}{2}-1)(\frac{n-1}{2}-2))^1,0^1\}, & \text{ if }~  2\nmid n \\
\{(-(\frac{n}{2}-3))^{(\frac{n}{2}-2)},((\frac{n}{2}-2)(\frac{n}{2}-3))^1,0^2\}, & \text{ if } ~  2\mid n
\end{cases}
$$
and
$$
\Ecn(\Gamma(G))=\begin{cases}
2(\frac{n-1}{2}-1)(\frac{n-1}{2}-2), & \text{ if }~  2\nmid n\\
2(\frac{n}{2}-2)(\frac{n}{2}-3), & \text{ if }~  2\mid n.
\end{cases}
$$
Hence, the result follows after simplification.
\end{proof}

\begin{thm}\label{Thm 3.2}
The CN-spectrum and CN-energy of commuting conjugacy class graph of the dicyclic group $T_{4n}=\langle x,y:~x^{2n}=1,x^n=y^2,y^{-1}xy=x^{-1} \rangle$ are given by
$$
\cnspec(\Gamma(T_{4n}))=\{(-(n-3))^{(n-2)},((n-2)(n-3))^1,0^2\}
$$
and
$$
\Ecn(\Gamma(T_{4n}))=2(n-2)(n-3).
$$
\end{thm}
\begin{proof}
By \cite[Proposition 2.2]{SA-CA-2020}, we have
$$
\Gamma(T_{4n})=\begin{cases}
K_{n-1}\cup 2K_1, & \text{ if } 2\mid n \\
K_{n-1}\cup K_2, & \text{ if } 2\nmid n. 
\end{cases}
$$
Applying Theorem \ref{CNS2}, we get
$$\cnspec(\Gamma(T_{4n}))=\{(-(n-3))^{(n-2)},((n-2)(n-3))^1,0^2\}$$
and
$$\Ecn(\Gamma(T_{4n}))=2(n-2)(n-3).$$
\end{proof}

\begin{thm}\label{Thm 3.5}
	The CN-spectrum and CN-energy of commuting conjugacy class graph of the semidihedral group $SD_{8n}=\langle x,y:~x^{4n}=y^2=1,yxy=x^{2n-1}\rangle$ are given by
	$$\cnspec(\Gamma(SD_{8n}))=\begin{cases}
		\{(-(2n-3))^{(2n-2)},((2n-2)(2n-3))^1,0^2\},& \text{if}~2\mid n \\
		\{(-(2n-4))^{(2n-3)},((2n-3)(2n-4))^1,-2^3,6^1\},& \text{if}~2\nmid n
	\end{cases}$$
	and
	$$\Ecn(\Gamma(SD_{8n}))=\begin{cases}
		2(2n-2)(2n-3),& \text{if}~2\mid n\\
		2(2n-3)(2n-4)+12,& \text{if}~2\nmid n.
	\end{cases}$$
\end{thm}
\begin{proof}
	By \cite[Proposition 2.5]{SA-CA-2020}, we have
	$$\Gamma(SD_{8n})=\begin{cases}
		K_{2n-1}\cup 2K_1,& \text{if}~2\mid n\\
		K_{2n-2}\cup K_4,& \text{if}~2\nmid n. 
	\end{cases}$$
	Now applying Theorem \ref{CNS2}, we get
	\begin{align*}
		&\cnspec(\Gamma(SD_{8n}))=\\
		&\quad \quad \begin{cases}
			\{(-((2n-1)-2))^{((2n-1)-1)},(((2n-1)-1)((2n-1)-2))^1,0^2\},& \text{if}~2\mid n\\
			\{(-((2n-2)-2))^{((2n-2)-1)},(((2n-2)-1)((2n-2)-2))^1,\\
			\quad \quad\quad \quad\quad \quad\quad \quad\quad \quad\quad \quad\quad \quad(-(4-2))^{(4-1)},((4-1)(4-2))^1\},& \text{if}~2\nmid n
		\end{cases}
	\end{align*}
	and
	$$
	\Ecn(\Gamma(V_{8n}))=\begin{cases}
		2((2n-1)-1)((2n-1)-2)+2\times2(1-1)(1-2),& \text{if}~2\mid n\\
		2((2n-2)-1)((2n-2)-2)+2(4-1)(4-2),& \text{if}~2\nmid n.
	\end{cases}
	$$
	Hence,  we get the required result on simplification.
\end{proof}

\begin{thm}\label{Thm 3.3}
The CN-spectrum and CN-energy of commuting conjugacy class graph of the group $U_{(n,m)}=\langle x,y:~x^{2n}=y^m=1,~x^{-1}yx=y^{-1}\rangle$ are given by
\begin{align*}
&\cnspec(\Gamma(U_{(n,m)})) =\\
&\quad\begin{cases}
	\{(-(n-2))^{2(n-1)},(n^2-3n+2)^2,(-(\frac{1}{2}(mn-2n-4))^{\frac{1}{2}(mn-2n-2)},\\ 
\quad \quad\quad \quad\quad \quad\quad \quad \quad \quad\quad \quad\quad \quad\quad	(\frac{1}{4}(mn-2n-2)(mn-2n-4))^1\},\text{ if }~  2\mid m\\
\{(-(n-2))^{(n-1)},(n^2-3n+2)^1,(-\frac{1}{2}(mn-n-4))^{\frac{1}{2}(mn-n-2)},\\ \quad \quad\quad \quad\quad \quad\quad \quad \quad \quad\quad \quad\quad \quad\quad(\frac{1}{4}(mn-n-2)(mn-n-4))^1\}, \quad \text{ if }~  2\nmid m
\end{cases}
\end{align*}
and
$$
\Ecn(\Gamma(U_{(n,m)}))=\begin{cases}
4(n^2-3n+2)+\frac{1}{2}(mn-2n-2)(mn-2n-4),~~~\text{ if }~  2\mid m\\
2(n^2-3n+2)+\frac{1}{2}(mn-n-2)(mn-n-4),~\quad \text{ if }~  2\nmid m.
\end{cases}
$$

\end{thm}
\begin{proof}
By \cite[Proposition 2.3]{SA-CA-2020}, we have
$$\Gamma(U_{(n,m)})=\begin{cases}
2K_n\cup K_{n(\frac{m}{2}-1)}, & \text{ if }~  2\mid m \\
K_n\cup K_{n(\frac{m-1}{2})}, & \text{ if }~  2\nmid m. 
 \end{cases}$$
Now applying Theorem \ref{CNS2}, we get
\begin{align*}
&\cnspec(\Gamma(U_{(n,m)}))=\\
&\quad \begin{cases}
\{(-(n-2))^{2(n-1)},((n-1)(n-2))^2,(-(n(\frac{m}{2}-1)-2))^{(n(\frac{m}{2}-1)-1)},\\ \quad \quad\quad \quad\quad \quad\quad \quad \quad \quad\quad \quad\quad \quad\quad((n(\frac{m}{2}-1)-1)(n(\frac{m}{2}-1)-2))\},~~~~~~~~~ \text{ if }~  2\mid m\\
\{(-(n-2))^{(n-1)},((n-1)(n-2))^1,(-(n(\frac{m-1}{2})-2))^{(n(\frac{m-1}{2})-1)},\\
\quad \quad\quad \quad\quad \quad\quad \quad \quad \quad\quad \quad\quad \quad\quad\quad ((n(\frac{m-1}{2})-1)(n(\frac{m-1}{2})-2)\}, ~\text{ if }~  2\nmid m.
\end{cases}
\end{align*}
and
$$
\Ecn(\Gamma(U_{(n,m)}))=\begin{cases}
2\times 2(n-1)(n-2)+2(n(\frac{m}{2}-1)-1))(n(\frac{m}{2}-1)-2),  \text{ if }~  2\mid m\\
2(n-1)(n-2)+2(n(\frac{m-1}{2})-1))(n(\frac{m-1}{2})-2), \qquad\,\,\text{ if }~  2\nmid m.
\end{cases}
$$
Hence, the result follows on simplification.
\end{proof}

\begin{cor}\label{U_6n}	The CN-spectrum and CN-energy of commuting conjugacy class graph of the group $U_{6n} =\langle x,y:~x^{2n}=y^3=1,~x^{-1}yx=y^{-1}\rangle$ are given by	
\[
\cnspec(\Gamma(U_{6n}))=\{(-(n-2))^{2(n-1)}, ((n-1)(n-2))^2\} 
\]
and
$
\Ecn(\Gamma(U_{6n}))=4(n-1)(n-2).
$
\end{cor}
\begin{proof}
The result follows from Theorem \ref{Thm 3.3} noting that $U_{(n,3)} = U_{6n}$.  
\end{proof}

\begin{thm}\label{Thm 3.4}
The CN-spectrum and CN-energy of commuting conjugacy class graph of the group $V_{8n}=\langle x,y:~x^{2n}=y^4=1,~yx=x^{-1}y^{-1},~y^{-1}x=x^{-1}y \rangle$ are given by
$$\cnspec(\Gamma(V_{8n}))=\begin{cases}
\{(-(2n-4))^{(2n-3)},((2n-3)(2n-4))^1,0^4\},& \text{ if }2\mid n\\
\{(-(2n-3))^{(2n-2)},((2n-2)(2n-3))^1,0^2\},& \text{ if }2\nmid n
\end{cases}$$
 and
 $$\Ecn(\Gamma(V_{8n}))=\begin{cases}
 2(2n-3)(2n-4),& \text{ if }2\mid n\\
 2(2n-2)(2n-3),& \text{ if }2\nmid n.
 \end{cases}$$
\end{thm}
\begin{proof}
By \cite[Proposition 2.4]{SA-CA-2020}, we have
$$\Gamma(V_{8n})=\begin{cases}
K_{2n-2}\cup 2K_2,& \text{ if }2\mid n\\
K_{2n-1}\cup 2K_1,& \text{ if }2\nmid n. 
 \end{cases}$$
Now applying Theorem \ref{CNS2}, we get
\begin{align*}
&\cnspec(\Gamma(V_{8n}))=\\
&\quad \begin{cases}
\{(-((2n-2)-2))^{((2n-2)-1)},(((2n-2)-1)((2n-2)-2))^1,0^2,0^2\},& \text{ if }2\mid n\\
\{(-((2n-1)-2))^{((2n-1)-1)},(((2n-1)-1)((2n-1)-2))^1,0^2\},& \text{ if }2\nmid n
\end{cases}
\end{align*}
and
$$
\Ecn(\Gamma(V_{8n}))=\begin{cases}
 2((2n-2)-1)((2n-2)-2),& \text{ if }2\mid n\\
 2((2n-1)-1)((2n-1)-2),& \text{ if }2\nmid n.
 \end{cases}
$$
Hence, the result follows after simplification.
\end{proof}

\begin{thm}\label{Thm 3.6}
The CN-spectrum and CN-energy of commuting conjugacy class graph of the group $G(p,m,n)=\langle x,y:~x^{p^m}=y^{p^n}=[x,y]^p=1,[x,[x,y]]=[y,[x,y]]=1\rangle$ are given by
\begin{align*}
    \cnspec(\Gamma(G(p,m,n)))=&\biggl\{(-(p^{m+ n-1}-p^{m+ n-2}-2))^{2(p^{m+ n-1}-p^{m+ n-2}-1)},\\
                              &((p^{m + n - 1}-p^{m+ n-2}-1)(p^{m + n -1}-p^{m +n-2}-2))^2,\\
                              &(-(p^{m}-p^{m-1}-2))^{(p^n-p^{n-1})(p^m-p^{m-1}-1)},\\
                              &((p^m-p^{m-1}-1)(p^m-p^{m-1}-2))^{(p^n-p^{n-1})}\biggl\}
\end{align*}
and
\begin{align*}
    \Ecn(\Gamma(G(p,m,n)))=&4(p^{m +n -1}-p^{m+n-2}-1)(p^{m+n-1}-p^{m +n-2}-2)\\
                           & +2(p^n-p^{n-1})(p^m-p^{m-1}-1)(p^m-p^{m-1}-2).
\end{align*}

\end{thm}
\begin{proof}
By \cite[Proposition 2.6]{SA-CA-2020}, we have
$$\Gamma(G(p,m,n)) = K_{p^{m-1}(p^n-p^{n-1})}\cup K_{p^{n-1}(p^m-p^{m-1})}\cup (p^n-p^{n-1})K_{p^{m-n}(p^n-p^{n-1})}.$$
Now applying Theorem \ref{CNS2}, we get
\begin{align*}
    \cnspec(\Gamma(G(p,m,n)))=&\biggl\{(-(p^{m-1}(p^n-p^{n-1})-2))^{(p^{m-1}(p^n-p^{n-1})-1)},\\
                              &((p^{m-1}(p^n-p^{n-1})-1)(p^{m-1}(p^n-p^{n-1})-2))^1,\\
                              &(-(p^{n-1}(p^m-p^{m-1})-2))^{(p^{n-1}(p^m-p^{m-1})-1)},\\
                              &((p^{n-1}(p^m-p^{m-1})-1)(p^{n-1}(p^m-p^{m-1})-2))^1,\\
                              &(-(p^{m-n}(p^n-p^{n-1})-2))^{(p^n-p^{n-1})(p^{m-n}(p^n-p^{n-1})-1)},\\
                              &((p^{m-n}(p^n-p^{n-1})-1)(p^{m-n}(p^n-p^{n-1})-2))^{(p^n-p^{n-1})}\biggr\}
\end{align*}
\begin{align*}
\text{ and }    \Ecn(\Gamma(G(p,m,n)))=&2(p^{m-1}(p^n-p^{n-1})-1)(p^{m-1}(p^n-p^{n-1})-2)+\\
                           &2(p^{n-1}(p^m-p^{m-1})-1)(p^{n-1}(p^m-p^{m-1})-2)+\\
                           &2(p^n-p^{n-1})(p^{m-n}(p^n-p^{n-1})-1)(p^{m-n}(p^n-p^{n-1})-2).
\end{align*}
Hence, the result follows.
\end{proof}
\subsection{Certain groups with given central quotients}
In this subsection we mainly consider finite groups whose central quotients are isomorphic to a group of order $p^2$, $p^3$ or a dihedral group  of order $2n$  and compute CN-spectrum and CN-energy of their  commuting conjugacy class graphs. 
\begin{thm}\label{Thm 3.7}
Let $G$ be a non-abelian finite group with centre $Z=Z(G)$ and $\frac{G}{Z}\cong Z_p\times Z_p$, where $p$ is a prime. Then 
$$\cnspec(\Gamma(G))=\{(-(n-2))^{(p+1)(n-1)},((n-1)(n-2))^{(p+1)}\}$$
and
$$\Ecn(\Gamma(G))=2(p+1)(n-1)(n-2),$$ where $n=\frac{(p-1)|Z|}{p}$.
\end{thm}
\begin{proof}
By \cite[Theorem 3.1]{SA-2020}, we have
$$\Gamma (G)=(p+1)K_n,~~~ \text{ where } n=\frac{(p-1)|Z|}{p}.$$
Hence, the result follows from Theorem \ref{CNS2}.
\end{proof}

\begin{cor}\label{Thm 3.8}
If $G$ is a non-abelian $p$-group of order $p^n$ and $|Z(G)|=p^{n-2}$, $p$ is prime and $n\geq 3$, then
\begin{align*}
&\cnspec(\Gamma(G))=\\
&\quad\{(-(p^{n-2}-p^{n-3}-2))^{(p+1)(p^{n-2}-p^{n-3}-1)},((p^{n-2}-p^{n-3}-1)(p^{n-2}-p^{n-3}-2))^{(p+1)}\}
\end{align*}
and 
\[
\Ecn(\Gamma(G))=2(p+1)(p^{n-2}-p^{n-3}-1)(p^{n-2}-p^{n-3}-2).
\]
\end{cor}
\begin{proof}
We have $\frac{G}{Z}\cong Z_p\times Z_p$ and $ n=\frac{(p-1)|Z|}{p} = (p-1)p^{n - 3}$. Hence, the result follows from Theorem \ref{Thm 3.7}.
\end{proof}

\begin{thm}\label{Thm 3.9}
Let $G$ be a non-abelian group with centre $Z$ such that $|\frac{G}{Z}|=p^3$, for a prime $p$. Then one of the following is satisfied:
\begin{enumerate}
    \item If $\frac{G}{Z}$ is abelian then  
    \begin{enumerate}
    	\item $\cnspec(\Gamma(G)) \, = \, \{(-(m-2))^{(m-1)}, \, ((m-1)(m-2))^1, \, (-(n-2))^{p^2(n-1)}, \,$ $((n-1)(n-2))^{p^2}\}$ and 
    	   	$\Ecn(\Gamma(G))=2(m-1)(m-2)+2p^2(n-1)(n-2)$ if $\Gamma(G)=K_m\cup p^2K_n$.
    
 \vspace{.3cm} 
     
    	\item  $\cnspec(\Gamma(G))= \{(-(n-2))^{(p^2+p+1)(n-1)},((n-1)(n-2))^{(p^2+p+1)}\}$ and   $\Ecn(\Gamma(G))= 2(p^2+p+1)(n-1)(n-2)$ if $\Gamma(G)=(p^2+p+1)K_n$.
    \end{enumerate}
 Here $m=\frac{(p^2-1)|Z|}{p}$ and $n=\frac{(p-1)|Z|}{p^2}$.
    \item If $\frac{G}{Z}$ is non-abelian then 
    \begin{enumerate}
    	\item $\cnspec(\Gamma(G))=\{(-(m-2))^{(m-1)},((m-1)(m-2))^1,(-(n_1-2))^{kp(n_1-1)},$  $((n_1-1)(n_1-2))^{kp},(-(n_2-2))^{(p-k)(n_2-1)},((n_2-1)(n_2-2))^{(p-k)}\}$ and
    	
    	 $\Ecn(\Gamma(G))=2(m-1)(m-2)+2kp(n_1-1)(n_1-2)+2(p-k)(n_2-1)(n_2-2)$ if  $\Gamma(G)=K_m\cup kpK_{n_1}\cup (p-k)K_{n_2}$. 

 \vspace{.3cm}   

   \item $\cnspec(\Gamma(G)) \, = \, \{(-(n_1-2))^{(kp+1)(n_1-1)}, \, ((n_1-1)(n_1-2))^{(kp+1)}, \quad $ $(-(n_2-2))^{(p+1-k)(n_2-1)},((n_2-1)(n_2-2))^{(p+1-k)}\}$ and 
   
   $\Ecn(\Gamma(G))=  2(kp+1)(n_1-1)(n_1-2)+2(p+1-k)(n_2-1)(n_2-2)$ if  $\Gamma(G)= (kp + 1)K_{n_1}\cup (p + 1 - k)K_{n_2}$. 

 \vspace{.3cm}     
  
     \item $\cnspec(\Gamma(G))=\{(-(m-2))^{(m-1)},((m-1)(m-2))^1,(-(n_2-2))^{p(n_2-1)},$ $ ((n_2-1)(n_2-2))^p\}$ and
           $\Ecn(\Gamma(G))=2(m-1)(m-2)+2p(n_2-1)(n_2-2)$ if  $\Gamma(G)= K_m\cup pK_{n_2}$. 

\vspace{.3cm}       

 \item $\cnspec(\Gamma(G))=\{(-(n_1-2))^{(p^2+p+1)(n_1-1)},((n_1-1)(n_1-2))^{(p^2+p+1)}\}$ and $\Ecn(\Gamma(G))=2(p^2+p+1)(n_1-1)(n_1-2)$ if  $\Gamma(G)= (p^2 + p + 1)K_{n_1}$.        

\vspace{.3cm}  
 
       \item $\cnspec(\Gamma(G)) \,\,= \,\, \{(-(n_1-2))^{(n_1-1)}, \,\, ((n_1-1)(n_1-2))^1,$\\ \indent \qquad\qquad\qquad\qquad\qquad\qquad $(-(n_2-2))^{(p+1)(n_2-1)},$ $((n_2-1)(n_2-2))^{(p+1)}\}$ and
        $\Ecn(\Gamma(G))=2(n_1-1)(n_1-2)+2(p+1)(n_2-1)(n_2-2)$ if $\Gamma(G)= K_{n_1}\cup (p + 1)K_{n_2}$. 
    \end{enumerate}
Here $m=\frac{(p^2-1)|Z|}{p}$, $n_1=\frac{(p-1)|Z|}{p^2}$ and $n_2=\frac{(p-1)|Z|}{p}$, $1\leq k\leq p$.
\end{enumerate}
\end{thm}
\begin{proof}

(a) If $\frac{G}{Z}$ is abelian then, by \cite[Theorem 3.3]{SA-2020}, we have $\Gamma(G)=K_m\cup p^2K_n$ or $(p^2+p+1)K_n$, where $m=\frac{(p^2-1)|Z|}{p}$ and $n=\frac{(p-1)|Z|}{p^2}$. 

If $\Gamma(G)=K_m\cup p^2K_n$ then, using Theorem \ref{CNS2}, we get
$$
\cnspec(\Gamma(G))=\{(-(m-2))^{(m-1)},(m-1)(m-2),(-(n-2))^{p^2(n-1)},((n-1)(n-2))^{p^2}\}
$$ 
and 
$$
\Ecn(\Gamma(G))=2(m-1)(m-2)+2p^2(n-1)(n-2).
$$

If $\Gamma(G)= (p^2+p+1)K_n$  then, using Theorem \ref{CNS2}, we get
$$
\cnspec(\Gamma(G))= \{(-(n-2))^{(p^2+p+1)(n-1)},((n-1)(n-2))^{(p^2+p+1)}\}
$$ 
and   
$$
\Ecn(\Gamma(G))= 2(p^2+p+1)(n-1)(n-2).
$$
This completes the proof of part (a).

(b) 
 If $\frac{G}{Z}$ is non-abelian then, by \cite[Theorem 3.3]{SA-2020}, we have $\Gamma(G)=K_m\cup kpK_{n_1}\cup (p-k)K_{n_2}$, $(kp+1)K_{n_1}\cup (p+1-k)K_{n_2}$, $K_m\cup pK_{n_2}$, $(p^2+p+1)K_{n_1}$ or $K_{n_1}\cup (p+1)K_{n_2}$, where $m=\frac{(p^2-1)|Z|}{p}$,$n_1=\frac{(p-1)|Z|}{p^2}$, $n_2=\frac{(p-1)|Z|}{p}$, $1\leq k\leq p$.

If $\Gamma(G)=K_m\cup kpK_{n_1}\cup (p-k)K_{n_2}$ then,  using Theorem \ref{CNS2}, we get
\begin{align*}
\cnspec(\Gamma(G))= &\{(-(m-2))^{(m-1)},((m-1)(m-2))^1,(-(n_1-2))^{kp(n_1-1)},\\
&  ((n_1-1)(n_1-2))^{kp},(-(n_2-2))^{(p-k)(n_2-1)},((n_2-1)(n_2-2))^{(p-k)}\}
\end{align*} 
and
$$
\Ecn(\Gamma(G))=2(m-1)(m-2)+2kp(n_1-1)(n_1-2)+2(p-k)(n_2-1)(n_2-2).
$$

If $\Gamma(G)=(kp+1)K_{n_1}\cup (p+1-k)K_{n_2}$ then,  using Theorem \ref{CNS2}, we get
\begin{align*}
\cnspec(\Gamma(G))=&\{(-(n_1-2))^{(kp+1)(n_1-1)},((n_1-1)(n_1-2))^{(kp+1)},\\
&\quad\quad\quad(-(n_2-2))^{(p+1-k)(n_2-1)},((n_2-1)(n_2-2))^{(p+1-k)}\}
\end{align*} 
and 
$$
\Ecn(\Gamma(G))=  2(kp+1)(n_1-1)(n_1-2)+2(p+1-k)(n_2-1)(n_2-2).
$$

If $\Gamma(G)=K_m\cup pK_{n_2}$ then,  using Theorem \ref{CNS2}, we get
\begin{align*}
\cnspec(\Gamma(G))=\{(-(m-2))^{(m-1)},&((m-1)(m-2))^1,\\
&(-(n_2-2))^{p(n_2-1)}, ((n_2-1)(n_2-2))^p\}
\end{align*} 
and
$$
\Ecn(\Gamma(G))=2(m-1)(m-2)+2p(n_2-1)(n_2-2).
$$

If $\Gamma(G)=(p^2+p+1)K_{n_1}$  then,  using Theorem \ref{CNS2}, we get
$$
\cnspec(\Gamma(G))=\{(-(n_1-2))^{(p^2+p+1)(n_1-1)},((n_1-1)(n_1-2))^{(p^2+p+1)}\}
$$ 
and 
$$
\Ecn(\Gamma(G))=2(p^2+p+1)(n_1-1)(n_1-2).
$$

If $\Gamma(G)=K_{n_1}\cup (p+1)K_{n_2}$ then,  using Theorem \ref{CNS2}, we get
\begin{align*}
\cnspec(\Gamma(G))=\{(-(n_1-2))^{(n_1-1)},&((n_1-1)(n_1-2))^1,\\
&(-(n_2-2))^{(p+1)(n_2-1)},((n_2-1)(n_2-2))^{(p+1)}\}
\end{align*}
and
$$
\Ecn(\Gamma(G))=2(n_1-1)(n_1-2)+2(p+1)(n_2-1)(n_2-2).
$$
This completes the proof of part (b).
\end{proof}

\begin{cor}\label{Thm 3.10}
Let $G$ be a non-abelian $p$-group of order $p^n$ and $|Z|=p^{n-3}$, where $p$ is prime and $n\geq 4$. Then one of the following are satisfied:
\begin{enumerate}
\item If $\frac{G}{Z}$ is abelian then 
\begin{enumerate}
	\item  
       $\cnspec(\Gamma(G))= \{(-(p^{n-2}-p^{n-4}-2))^{(p^{n-2}-p^{n-4}-1)},((p^{n-2}-p^{n-4}-1)(p^{n-2}-p^{n-4}-2))^1, (-(p^{n-4}-p^{n-5}-2))^{p^2(p^{n-4}-p^{n-5}-1)},((p^{n-4}-p^{n-5}-1)(p^{n-4}-p^{n-5}-2))^{p^2}\}$ and 
       \begin{align*}
       	       \Ecn(\Gamma(G))=2(p^{n-2}&-p^{n-4}-1)(p^{n-2}-p^{n-4}-2) \\
       	       &+2p^2(p^{n-4}-p^{n-5}-1)(p^{n-4}-p^{n-5}-2)
       \end{align*}
   if $\Gamma(G)=K_{p^{n-4}(p^2 - 1)}\cup p^2K_{p^{n - 5}(p-1)}$.
\vspace{.3cm}   
      \item  $\cnspec(\Gamma(G))= \{(-(p^{n-4}-p^{n-5}-2))^{(p^2+p+1)(p^{n-4}-p^{n-5}-1)},
                          ((p^{n-4}-p^{n-5}-1)(p^{n-4}-p^{n-5}-2))^{(p^2+p+1)}\}$
   and
$$
\Ecn(\Gamma(G))=  2(p^2+p+1)(p^{n-4}-p^{n-5}-1)(p^{n-4}-p^{n-5}-2)
$$
if $\Gamma(G)=(p^2+p+1)K_{p^{n - 5}(p-1)}$.
\end{enumerate}
\vspace{.3cm}
\item If $\frac{G}{Z}$ is non-abelian then
\begin{enumerate}
	 \item    $\cnspec(\Gamma(G))=\{(-(p^{n-2}-p^{n-4}-2))^{(p^{n-2}-p^{n-4}-1)},((p^{n-2}-p^{n-4}-1)(p^{n-2}-p^{n-4}-2))^1,
	(-(p^{n-4}-p^{n-5}-2))^{kp(p^{n-4}-p^{n-5}-1)},((p^{n-4}-p^{n-5}-1)(p^{n-4}-p^{n-5}-2))^{kp},
	(-(p^{n-3}-p^{n-4}-2))^{(p-k)(p^{n-3}-p^{n-4}-1)},((p^{n-3}-p^{n-4}-1)(p^{n-3}-p^{n-4}-2))^{(p-k)}\}$ and
	
	$\Ecn(\Gamma(G))= 2(p^{n-2}-p^{n-4}-1)(p^{n-2}-p^{n-4}-2)+2kp(p^{n-4}-p^{n-5}-1)(p^{n-4}-p^{n-5}-2)+
	2(p-k)(p^{n-3}-p^{n-4}-1)(p^{n-3}-p^{n-4}-2)$
	if  $\Gamma(G)= K_{p^{n-4}(p^2-1)} \cup kpK_{p^{n-5}(p-1)}\cup (p-k)K_{p^{n-4}(p-1)}$.        
	\vspace{.3cm}
\vspace{.3cm}                       
\item    $\cnspec(\Gamma(G))=\{(-(p^{n-4}-p^{n-5}-2))^{(kp+1)(p^{n-4}-p^{n-5}-1)},((p^{n-4}-p^{n-5}-1)(p^{n-4}-p^{n-5}-2))^{(kp+1)},
(-(p^{n-3}-p^{n-4}-2))^{(p+1-k)(p^{n-3}-p^{n-4}-1)},((p^{n-3}-p^{n-4}-1)(p^{n-3}-p^{n-4}-2))^{(p+1-k)}\}$ and

$\Ecn(\Gamma(G))= 2(kp+1)(p^{n-4}-p^{n-5}-1)(p^{n-4}-p^{n-5}-2)+2(p+1-k)(p^{n-3}-p^{n-4}-1)(p^{n-3}-p^{n-4}-2)$ if $\Gamma(G)= (kp+1)K_{p^{n-5}(p-1)}\cup (p + 1 - k)K_{p^{n-4}(p-1)}$. 
 
\item    $\cnspec(\Gamma(G))=\{(-(p^{n-2}-p^{n-4}-2))^{(p^{n-2}-p^{n-4}-1)},((p^{n-2}-p^{n-4}-1)(p^{n-2}-p^{n-4}-2))^1,                     (-(p^{n-3}-p^{n-4}-2))^{p(p^{n-3}-p^{n-4}-1)},((p^{n-3}-p^{n-4}-1)(p^{n-3}-p^{n-4}-2))^p\}$ and

 $\Ecn(\Gamma(G))=2(p^{n-2}-p^{n-4}-1)(p^{n-2}-p^{n-4}-2)+2p(p^{n-3}-p^{n-4}-1)(p^{n-3}-p^{n-4}-2)$ if  $\Gamma(G)= K_{p^{n - 4}(p^2 - 1)}\cup pK_{p^{n-4}(p-1)}$. 
\vspace{.3cm}
\item    $\cnspec(\Gamma(G))=\{(-(p^{n-4}-p^{n-5}-2))^{(p^2+p+1)(p^{n-4}-p^{n-5}-1)},((p^{n-4}-p^{n-5}-1)(p^{n-4}-p^{n-5}-2))^{(p^2+p+1)}\}$ and

 $\Ecn(\Gamma(G))=2(p^2+p+1)(p^{n-4}-p^{n-5}-1)(p^{n-4}-p^{n-5}-2)$ if  $\Gamma(G)= (p^2 + p + 1)K_{p^{n-5}(p-1)}$. 
\vspace{.3cm}
\item    $\cnspec(\Gamma(G))=\{(-(p^{n-4}-p^{n-5}-2))^{(p^{n-4}-p^{n-5}-1)},((p^{n-4}-p^{n-5}-1)(p^{n-4}-p^{n-5}-2))^1, 
(-(p^{n-3}-p^{n-4}-2))^{(p+1)(p^{n-3}-p^{n-4}-1)},((p^{n-3}-p^{n-4}-1)(p^{n-3}-p^{n-4}-2))^{(p+1)}\}$ and

 $\Ecn(\Gamma(G))=2(p^{n-4}-p^{n-5}-1)(p^{n-4}-p^{n-5}-2)+2(p+1)(p^{n-3}-p^{n-4}-1)(p^{n-3}-p^{n-4}-2)$, if  $\Gamma(G)= K_{p^{n-5}(p-1)}\cup (p+1)K_{p^{n-4}(p-1)}$,
 \end{enumerate}
where $1\leq k\leq p$.
\end{enumerate}

\end{cor}
\begin{proof}
We have $|\frac{G}{Z}|=p^3$ and   $m=\frac{(p^2-1)|Z|}{p} = (p^2-1)p^{n-4}$, $n = n_1=\frac{(p-1)|Z|}{p^2} = (p-1)p^{n-5}$ and $n_2=\frac{(p-1)|Z|}{p} = (p-1)p^{n-4}$.
Hence, the result follows from Theorem \ref{Thm 3.9}.
%
\end{proof}

\begin{cor}\label{Thm 3.11}
Let $G$ be a non-abelian $p$-group of order $p^4$. Then 
\begin{enumerate}
 \item   $\cnspec(\Gamma(G))=\{(-(p^2-p-2))^{(p+1)(p^2-p-1)},((p^2-p-1)(p^2-p-2))^{(p+1)}\}$ and $\Ecn(\Gamma(G))=2(p+1)(p^2-p-1)(p^2-p-2)$ if $\Gamma(G) = (p+1)K_{p(p-1)}$.
     
\vspace{.3cm}     
\item   $\cnspec(\Gamma(G))=     \{(-(p^2-3))^{(p^2-2)},((p^2-2)(p^2-3))^1,(-(p-3))^{(p^2-2p)},((p-2)(p-3))^p\}$ and  $\Ecn(\Gamma(G))=2(p^2-2)(p^2-3)+2p(p-2)(p-3)$ if $\Gamma(G)= K_{(p^2-1)}\cup pK_{p-1}$.
\end{enumerate}
\end{cor}
\begin{proof}
If $G$ is a non-abelian $p$-group of order $p^4$ then $|Z(G)| = p$ or $p^2$. Suppose that $|Z(G)| = p^2$. Then by Corollary \ref{Thm 3.8}, we get 
$$
\cnspec(\Gamma(G))=\{(-(p^2-p-2))^{(p+1)(p^2-p-1)},((p^2-p-1)(p^2-p-2))^{(p+1)}\}
$$ and
$$
\Ecn(\Gamma(G))=2(p+1)(p^2-p-1)(p^2-p-2).
$$
In this case,  by  \cite[Corollary 3.2]{SA-2020}, it follows that $\Gamma(G) = (p+1)K_{p(p-1)}$.

If $|Z(G)| = p$ then, by  \cite[Corollary 3.5]{SA-2020}, we have $\Gamma(G)= K_{(p^2-1)}\cup pK_{p-1}$.   
Using Corollary \ref{Thm 3.10}(b), we get
\begin{align*}
    \cnspec(\Gamma(G))=\{(-(p^2-1-&2))^{(p^2-1-1)},((p^2-1-1)(p^2-1-2))^1,\\
    &(-(p-1-2))^{p(p-1-1)},((p-1-1)(p-1-2))^p\}
\end{align*}
and
\begin{align*}
    \Ecn(\Gamma(G))=& 2(p^2-1-1)(p^2-1-2)+2p(p-1-1)(p-1-2).
\end{align*}
After some simplification we get our required result.
\end{proof}


\begin{thm}\label{D2n-spec-energy}
Let $G$ be a non-abelian finite group with center $Z=Z(G)$ and $|Z|=z$. If   $\frac{G}{Z}\cong D_{2n}$. Then
\begin{align*}
	\cnspec&(\Gamma(G))=\\
	&\begin{cases}
		\left\lbrace\left(-\left(\frac{(n-1)z}{2}-2\right)\right)^{\left(\frac{(n-1)z}{2}-1\right)},\left(\left(\frac{(n-1)z}{2}\right)^2 - \frac{3(n-1)z}{2} + 2\right)^1,\right.\\
		\qquad\qquad\qquad 	\qquad\qquad	\left.\left(-\left(\frac{z}{2}-2\right)\right)^{\left(\frac{z}{2}-1\right)},\left(\left(\frac{z}{2}\right)^2 - \frac{3z}{2} + 2\right)^2\right\rbrace, &\text{ if } 2 \mid n\\		
		\left\lbrace\left(-\left(\frac{(n-1)z}{2}-2\right)\right)^{\left(\frac{(n-1)z}{2}-1\right)},\left(\left(\frac{(n-1)z}{2}\right)^2 - \frac{3(n-1)z}{2} + 2\right)^1,\right.\\
		\qquad\qquad\qquad 	\qquad\qquad	\left.(-(z-2))^{(z-1)},(z^2-3z+2)^1\right\rbrace,  &\text{ if } 2 \nmid n
	\end{cases}
\end{align*}	
	and 
\[
\Ecn(\Gamma(G))=\begin{cases}
		\frac{n^2 z^2}{2}-n z^2-3 n z+\frac{3 z^2}{2}-3 z+12, &\text{if $2 \mid n$}\\
		\frac{n^2 z^2}{2}-n z^2-3 n z+\frac{5 z^2}{2}-3 z+8, &\text{if $2 \nmid n$.}
	\end{cases}
\]
\end{thm}

\begin{proof}
	By \cite[Theorem 1.2]{Salah-2020}, we have
	$$\Gamma(G)=\begin{cases}
		K_{\frac{(n-1)z}{2}}\cup 2K_{\frac{z}{2}}, & \text{ if } 2 \mid n\\
		K_{\frac{(n-1)z}{2}}\cup K_{z}, & \text{ if } 2 \nmid n.
	\end{cases}$$
	Now applying Theorem \ref{CNS2}, we get
\begin{align*}
\cnspec&(\Gamma(G))=\\
&\begin{cases}
	\left\lbrace\left(-\left(\frac{(n-1)z}{2}-2\right)\right)^{\left(\frac{(n-1)z}{2}-1\right)},\left(\left(\frac{(n-1)z}{2}-2\right)\left(\frac{(n-1)z}{2}-1\right)\right)^1,\right.\\
		\qquad\qquad\qquad 	\qquad\qquad	\left.\left(-\left(\frac{z}{2}-2\right)\right)^{\left(\frac{z}{2}-1\right)},\left(\left(\frac{z}{2}-2\right) \left(\frac{z}{2}-1\right)\right)^2\right\rbrace, &\text{ if } 2 \mid n\\		
			\left\lbrace\left(-\left(\frac{(n-1)z}{2}-2\right)\right)^{\left(\frac{(n-1)z}{2}-1\right)},\left(\left(\frac{(n-1)z}{2}-2\right)\left(\frac{(n-1)z}{2}-1\right)\right)^1,\right.\\
	\qquad\qquad\qquad 	\qquad\qquad	(-(z-2))^{(z-1)},((z-2)(z-1))^1\biggr\},  &\text{ if } 2 \nmid n
	\end{cases}
\end{align*}
	and
	$$\Ecn(\Gamma(G))=\begin{cases}
		2(\frac{(n-1)z}{2}-1)(\frac{(n-1)z}{2}-2)+2\times 2(\frac{z}{2}-1)(\frac{z}{2}-2), &\text{ if } 2 \mid n\\
		2(\frac{(n-1)z}{2}-1)(\frac{(n-1)z}{2}-2)+2(z-1)(z-2), &\text{ if } 2 \nmid n.
	\end{cases}$$
Hence the result follows on simplification.
\end{proof}
\noindent We conclude this section with the following remark.
\begin{rem}\label{Remark}
	We have $\frac{D_{2n}}{Z(D_{2n})}\cong 	D_{n}$ or $D_{2n}$ according as $n$ is even or odd,   $\frac{T_{4n}}{Z(T_{4n})}\cong D_{2n}$ and $\frac{U_{6n}}{Z(U_{6n})}=D_{2\times 3}$. Therefore, Theorem \ref{Thm 3.1}, Theorem \ref{Thm 3.2} and Corollary \ref{U_6n} can also be obtained from Theorem \ref{D2n-spec-energy}.  
\end{rem}

\section{Some consequences} 
We begin this section with the following consequence of the results obtained in Section 2.
\begin{thm}
	\begin{enumerate}
		\item If $G$ is isomorphic to $D_{2n}$, $T_{4n}$, $SD_{8n}$, $U_{(n,m)}$, $U_{6n}$, $V_{8n}$ or $G(p, m, n)$ then  commuting conjugacy class graph of $G$ is CN-integral.
		\item If $G$ is a finite group and $\frac{G}{Z(G)}$ is isomorphic to $Z_p\times Z_p$ or $D_{2n}$ then  commuting conjugacy class graph of $G$  is CN-integral.
		\item If $G$ is a finite group and $\frac{G}{Z(G)}$ is of order $p^3$ then  commuting conjugacy class graph of $G$  is CN-integral.
\end{enumerate}	
\end{thm}

\begin{thm}
	Let $G$ be the semidihedral group $SD_{8n}$, where $n\geq 2$. Then the commuting conjugacy class graph of $G$ is not CN-hyperenergetic.
\end{thm}
\begin{proof}
We have $|V(\Gamma (G))| = 2n+1$ or $2n+2$ according as $n$ is even or odd. By Theorem \ref{CNS2} and Theorem \ref{Thm 3.5}, we get
\[
\Ecn(K_{|V(\Gamma (G))|})-\Ecn(\Gamma (G))= \begin{cases}
4n(2n-1) - 2(2n-3)(2n-2), &\text{ if } 2 \mid n\\
4n(2n+1) - 2(2n-4)(2n-3)-12, &\text{ if } 2 \nmid n.
\end{cases}
\]
Since  $n\geq 2$ we have
\[
4n(2n-1) - 2(2n-3)(2n-2) = 2(8n-6)>0
\]
and 
\[
4n(2n+1) - 2(2n-4)(2n-3)-12 = 2(16n-18)>0.
\]	
Therefore, $\Ecn(K_{|V(\Gamma (G))|}) > \Ecn(\Gamma (G))$ and so  $\Gamma (G)$ is not CN-hyperenergetic.
\end{proof}

\begin{thm}
Let $G$ be the group $U_{(n,m)}$, where $n\geq 2,~m\geq 3$. Then the commuting conjugacy class graph of $G$  is not CN-hyperenergetic.
\end{thm}
\begin{proof}
We have  $|V(\Gamma (G))|=n+\frac{nm}{2}$ or $\frac{n}{2}+\frac{nm}{2}$	according as $m$ is even or odd. If $m$ is even then, by Theorem \ref{CNS2} and Theorem \ref{Thm 3.3}, we have 
\begin{align*}
    \Ecn(K_{|V(\Gamma (G))|})-\Ecn(\Gamma (G))=& \frac{1}{2}(2n+ nm -2)(2n+nm -4) - 4(n^2-3n+2)\\
    & \qquad\qquad\qquad\qquad -\frac{1}{2}(mn-2n-2)(mn-2n-4)\\
    =&4(n^2(m-1)-2)>0,
\end{align*}
since $n\geq 2$ and $m\geq 3$.

If $m$ is odd then, by Theorem \ref{CNS2} and Theorem \ref{Thm 3.3}, we have 
\begin{align*}
    \Ecn(K_{|V(\Gamma (G))|})-\Ecn(\Gamma (G))=& \frac{1}{2}(n +nm - 2)(n +nm - 4) - 2(n^2-3n+2)\\
    &  \qquad\qquad\qquad\qquad - \frac{1}{2}(mn-n-2)(mn-n-4)\\
    =&2n^2(m - 1)-4 \geq 0,
\end{align*}
since $n\geq 2$ and $m\geq 3$. Hence, the result follows.
\end{proof}

\begin{thm}
Let $G$ be the group $V_{8n}$, where $n\geq 2$. Then the commuting conjugacy class graph of $G$ is not CN-hyperenergetic.
\end{thm}
\begin{proof}
We have  $|V(\Gamma (G))| = 2n+2$ or $2n+1$ according as $n$ is even or odd.
By  Theorem \ref{CNS2} and Theorem \ref{Thm 3.4}, we get
\[
\Ecn(K_{|V(\Gamma (G))|})-\Ecn(\Gamma (G))= \begin{cases}
	4n(2n+1)-2(2n-3)(2n-4), & \text{ if } 2 \mid n\\
	4n(2n-1)-2(2n-3)(2n-2), & \text{ if } 2 \nmid n.
\end{cases}
\]
Since $n \geq 2$ we have
\[
4n(2n+1)-2(2n-3)(2n-4) = 8(4n-3)>0
\]
and
\[
4n(2n-1)-2(2n-3)(2n-2) = 4(4n-3)>0.
\]
Therefore, $\Ecn(K_{|V(\Gamma (G))|}) > \Ecn(\Gamma (G))$ and so  $\Gamma (G)$ is not CN-hyperenergetic.
\end{proof}

\begin{thm}
The commuting conjugacy class graph of the group $G(p,m,n)$, where $p$ is prime, $m\geq 1,~n\geq 1$ is not CN-hyperenergetic.
\end{thm}
\begin{proof}
We have $|V(\Gamma (G(p,m,n)))|= p^{m+n}-p^{m+n-2}$.
Now by using Theorem \ref{CNS2} and Theorem \ref{Thm 3.6}, we get
\begin{align*}
    \Ecn&(K_{|V(\Gamma (G(p,m,n)))|}-\Ecn(\Gamma (G(p,m,n)))\\
    &\qquad\qquad= 2 p^{2 m+n-3}-6 p^{2 m+n-2}+6 p^{2 m+n-1}-2 p^{2 m+n}-2 p^{2 m+2 n-4}\\
    &\qquad\qquad\qquad + 8 p^{2 m+2 n-3}-8 p^{2 m+2 n-2}+2 p^{2 m+2 n}+4 p^{n-1}-4 p^n-4:= \beta.
\end{align*}
Let $S_1=8p^{2m+2n-3}+2p^{2m+2n}-2p^{2m+2n-4}-8p^{2m+2n-2}$ and  $S_2=6p^{2m+n-1}-6p^{2m+n-2}-4$. Then  
\[
S_1 = 8p^{2m+2n-4}(2p^2(p^2-4)+(8p-2)) >0,
\]
and 
\[
S_2=6p^{2m+n-1}(1-p^{-1})-4=6p^{2m+n-2}(p-1)-4 \geq 8,
\]
since   $p^2-4\geq 0$,  $8p-2\geq 14$ and $6p^{2m+n-2}(p-1)\geq 12$.

Let $S_3=2p^{2m+n-3}+4p^{n-1}-4p^n=2p^{n-3}(p^{2m}+2p^2-2p^3)$. If $m\geq 2$ then $2m-3\geq 1\implies p^{2m-3}-2\geq 0$. So, $p^{2m}+2p^2-2p^3=p^3(p^{2m-3}-2)+2p^2\geq 0$. Hence $S_3\geq 0$.

Let $S_4=S_1-2p^{2m+n}=p^{2m+n}(8p^{n-3}+2p^n-2p^{n-4}-8p^{n-2}-2)$. Then  
\[
8p^{n-3}+2p^n-2p^{n-4}-8p^{n-2}-2=2p^{n-4}(p^2(p^2-4)+4p-1)-2 \geq 14,
\]
if $n\geq 4$ (in this case   $2p^{n-4}\geq 2$ and $p^2(p^2-4)+4p-1\geq 7$). Therefore $S_4>0$.
Hence, $\beta=S_2+S_3+S_4>0$, for all $m\geq 2$ and $n\geq 4$.

We shall  now consider the following cases:


\noindent {\bf Case 1.} $m\geq 2$ and $n=1$
\par In this case $\beta=2p^{2m-1}(p-1)^2(p+1)-4p>0$, since $p+1\geq 3$ and $2m-1\geq 1$.

\noindent {\bf Case 2.} $m\geq 2$ and $~n=2$
\par In this case  $\beta=2p^{2m-1}(p^3(p^2-5)+p(7p-4)+1)-4(p^2-p+1)$. We have $2p^{2m-1}\geq 16$ since $m\geq 2$.
If  $p=2$ then  $\beta=13.2^{2m}-12>0$. Suppose that  $p\geq 3$. Then  $p^2-5\geq 4$, $7p-4\geq 17$ and $p^3\geq p^2$. Hence, $\beta > 0$. 

\noindent {\bf Case 3.} $m\geq 2$ and $n=3$
\par In this case  $\beta=2p^{2m}(1+3p(p^2-1)+2p^2)+2p^{2m+4}(p^2-4) -4p^3 + 4(p^2-1)$.  If  $p=2$ then  $\beta=27.2^{2m+1}-20>0$. 
If $p\geq 3$ then $p^4(p^2-4)>p^3$ and $2p^{2m} \geq 4$. Therefore,
$2p^{2m + 4}(p^2-4)>4p^3$.
Hence $\beta >0$. 

\noindent {\bf Case 4.} $m=1$ and $n\geq 1$
\par  In this case
\begin{align*}
    \beta=& 6p^{n-1}-4+2p^{n+1}(3-5p^{-1})+p^{n+2}(-2-2p^{n-4}+8p^{n-3}-8p^{n-2}+2p^n)\\
         =& 6p^{n-1}-4+2p^{n+1}(3-5p^{-1})+p^{n+2}(2p^{n-4}(p^2(p^2-4)+4p-1)-2).
\end{align*}
Note that $\beta > 0$ if  $n\geq 4$. If $n=1$ then $\beta= 2p^4-2p^3 -2p^2-2p = 2p(p(p(p-1)-1)-1)>0$, since $p -1 \geq 1$.
If $n=2$ then $\beta=6p-12p^2+14p^3-10p^4+2p^6-4=2p^3(7-6p^{-1})+2p^6(1+3p^{-5}-5p^{-2})-4$. 
It is clear that $\beta>0$ if $p\geq 3$. If $p=2$ then $\beta=40>0$.
If $n=3$ then $\beta=(6p^2-4)+2p^4(2-5p^{-1}+3p)+2p^8(1-4p^{-2})>0$.


Hence, $\Ecn(K_{|V(\Gamma (G))|}) > \Ecn(\Gamma (G))$ and so  $\Gamma (G)$ is not CN-hyperenergetic.
\end{proof}

\begin{thm}\label{D2n Theorem}
	Let $G$ be a finite group with centre $Z$ and $|Z|=z$ such that $\frac{G}{Z}\cong D_{2n}$. Then the commuting conjugacy class graph of $G$ is not CN-hyperenergetic. In particular, if $G\cong D_{6}$ then $\Gamma (G)$ is CN-borderenergetic.
\end{thm}
\begin{proof}
If $n$ is even then $|V(\Gamma(G))|=\frac{(n-1)z}{2}+2 \frac{z}{2}=\frac{n z}{2}+\frac{z}{2}$. By  Theorem \ref{CNS2} and Theorem \ref{D2n-spec-energy}, we get
	\begin{align*}
		\Ecn&(K_{|V(\Gamma (G))|})-\Ecn(\Gamma (G))\\
		& = 2 \left(\frac{n z}{2} + \frac{z}{2}-1\right) \left(\frac{n z}{2} + \frac{z}{2} - 2\right)
		- \left\lbrace\frac{n^2 z^2}{2} - n z^2 - 3 n z + \frac{3 z^2}{2} - 3 z + 12\right\rbrace\\
		& = z^2(2n - 1) - 8 > 0,~~~~\text{since}~n\geq 3 ~\text{and}~z\geq 2.
	\end{align*}
If $n$ is odd then  $|V(\Gamma(G))|=\frac{(n-1)z}{2}+z=\frac{n z}{2}+\frac{z}{2}$. By  Theorem \ref{CNS2} and Theorem \ref{D2n-spec-energy}, we get
	\begin{align*}
		\Ecn&(K_{|V(\Gamma (G))|})-\Ecn(\Gamma (G))\\
		& = 2 \left(\frac{n z}{2} + \frac{z}{2} - 1\right) \left(\frac{n z}{2} + \frac{z}{2} - 2\right)
		-  \left\lbrace\frac{n^2 z^2}{2} - n z^2 - 3 n z + \frac{5 z^2}{2} - 3 z + 8\right\rbrace\\
		& = 2 z^2(n - 1) - 4.
	\end{align*}
We have 
\[
2 z^2(n - 1) - 4
\begin{cases}
	=0,~~~~\text{ for } n = 3, z = 1\\
	>0,~~~~\text{ otherwise. }
\end{cases}
\]
Hence the result follows.
\end{proof}
In view of Theorem \ref{D2n Theorem} and  Remark \ref{Remark} we get the following corollary.
\begin{cor}
	The commuting conjugacy class graph of
	\begin{enumerate}
		\item the dihedral group $D_{2m}$, where $m\geq 3$, is not CN-hyperenergetic.
		\item the generalized quarternion group $T_{4n}$, where $n\geq 2$, is not CN-hyperenergetic.
		\item the group $U_{6n}$ is not CN-hyperenergetic.
	\end{enumerate}
\end{cor}

\begin{thm}\label{Realted to 3.7}
Let $G$ be a non-abelian finite group with centre $Z(G)$ and $\frac{G}{Z(G)}\cong Z_p\times Z_p$, where $p$ is a prime. Then the commuting conjugacy class graph of  $G$ is not CN-hyperenergetic. 
\end{thm}
\begin{proof}
We have  $|V(\Gamma (G))|=np+n$, where $n=\frac{(p-1)|Z(G)|}{p}$. By  Theorem \ref{CNS2} and Theorem \ref{Thm 3.7}, we get
\begin{align*}
    \Ecn(K_{|V(\Gamma (G))|})-\Ecn(\Gamma (G))=& 2 (np+n-1) (n p+n-2)-2 (n-1) (n-2) (p+1)\\
    =& 2n^2p+2p(n^2p-2) > 0,
\end{align*}
since  $p\geq 2$ and $n$ is a positive integer.
Hence, $\Gamma (G)$ is not CN-hyperenergetic.
\end{proof}

\begin{cor}
Let $G$ be a non-abelian $p$-group of order $p^n$ and $|Z(G)|=p^{n-2}$, where $p$ is a prime and $n\geq 3$. Then the commuting conjugacy class graph of  $G$ is not CN-hyperenergetic.
\end{cor}
\begin{proof}
If $G$ is a non-abelian $p$-group of order $p^n$ and $|Z(G)|=p^{n-2}$ then  $\frac{G}{Z(G)}\cong Z_p\times Z_p$,  $p$ is prime and $n\geq 3$. Hence, the result follows from Theorem \ref{Realted to 3.7}.
\end{proof}

\begin{thm} 
Let $G$ be a non-abelian $p$-group of order $p^n$ and $|Z(G)|=p^{n-3}$, where $p$ is prime and $n\geq 4$ then
\begin{enumerate}
   \item if $\frac{G}{Z(G)}$ is abelian then $\Gamma(G)$ is not CN-hyperenergetic,
    \item if $\frac{G}{Z(G)}$ is non-abelian then $\Gamma(G)$ is not CN-hyperenergetic.
\end{enumerate}
\end{thm}
\begin{proof}
(a) If  $\frac{G}{Z(G)}$ is abelian then, by Theorem \ref{CNS2} and Corollary \ref{Thm 3.10}, we get
\begin{align*}
\Ecn&(K_{|V(\Gamma (G))|})-\Ecn(\Gamma (G)) \\
&= \begin{cases}
	-4p^2+2p^{2n-6}(p^2-2)+2p^{2n-8}(2p^3(p-2)+4p-1) := \beta_1(p, n),\\
	\qquad\qquad\qquad\qquad\qquad\qquad\qquad \text{ if } \Gamma(G) = K_{p^{n-4}(p^2-1)} \cup p^2 K_{p^{n-5}(p-1)}\\
	  2p^{2n-7}(p^3-p-1) - (4p^2+4p) + 2p^{2n-9} := \beta_2(p, n), \\
	 \qquad\qquad\qquad\qquad\qquad\qquad\qquad  \text{ if } \Gamma(G) = (p^2+p+1)K_{p^{n-5}(p-1)}
\end{cases}
\end{align*}
\noindent noting that $|V( K_{p^{n-4}(p^2-1)} \cup p^2 K_{p^{n-5}(p-1)})| = 2 p^{n-2}  -p^{n-3}  -p^{n-4}$ and $|V((p^2+p+1)K_{p^{n-5}(p-1)})| = p^{n-2}-p^{n-5}$.

\vspace{.3cm}

\noindent \textbf{Case 1.} $\Ecn(K_{|V(\Gamma (G))|})-\Ecn(\Gamma (G)) = \beta_1(p, n)$

In this case $p \geq 2$ and $n \geq 5$. Since $2p^3(p-2) + 4p-1>0$ and $2p^{2n-6}(p^2-2)\geq 4p^2$ we have $\beta_1(p, n) > 0$. 
\vspace{.3cm}

\noindent \textbf{Case 2.} $\Ecn(K_{|V(\Gamma (G))|})-\Ecn(\Gamma (G)) = \beta_2(p, n)$

In this case $p \geq 2$ and $n \geq 5$. Since $p^{2n-7}\geq p^2$  and $p^3-p-1>4$ we have 
$p^{2n-7}(p^3-p-1)>4p^2$ and so $p^{2n-7}(p^3-p-1)>4p$. Therefore, $2p^{2n-7}(p^3-p-1)>4p+4p^2$. Hence, $\beta_2(p, n)>0$.

(b) If $\frac{G}{Z}$ is non-abelian then, by Theorem \ref{CNS2} and Corollary \ref{Thm 3.10}, we get
 \begin{align*}
\Ecn&(K_{|V(\Gamma (G))|})-\Ecn(\Gamma (G))\\
&=\begin{cases}
2kp^{2n-8}(3-p^{-1})+2kp{2n-6}(1-3p^{-1})+2p^{2n-4}(3-5p^{-1}) +4k\\
\qquad +2p(p^{2n-8}-2)+2p(p^{2n-7}-2k):=\mu_1(p, n, k),\\
 \qquad\qquad\qquad\qquad	\text{ if }  \Gamma(G)= K_{p^{n-4}(p^2-1)} \cup kpK_{p^{n-5}(p-1)}\cup (p-k)K_{p^{n-4}(p-1)}\vspace{.3cm}\\
 4(k-1)+\{2p^{2n-4}(1-p^{-1})-4p^{2n-8}\}+\{2p^{2n-6}(1-p^{-1})-4p\}\\
\, +\{2kp^{2n-6}(1-3p^{-1})-4kp\}+2kp^{2n-8}(3-p^{-1})+4p^{2n-9}:=\mu_2(p, n. k),\\
 \qquad\qquad\qquad\qquad \text{ if } \Gamma(G)= (kp+1)K_{p^{n-5}(p-1)}\cup (p + 1 - k)K_{p^{n-4}(p-1)}\vspace{.3cm}\\
	2p^{2 n-4}(3-5p^{-1})+2p(p^{2n-6}+p^{2n-5}-2) := \mu_3(p, n), \\
\qquad\qquad\qquad\qquad	\text{ if } \Gamma(G)=K_{p^{n-4}(p^2-1)}\cup pK_{p^{n-4}(p-1)} \vspace{.3cm}\\
 2p^{2n-7}(p^3-p-1) - (4p^2+4p) + 2p^{2n-9} :=\mu_4(p, n),\\ 
\qquad\qquad\qquad\qquad \text{ if }  \Gamma(G)= (p^2+p+1)K_{p^{n-5}(p-1)}\vspace{.3cm}\\
4p^{2n-9} - 4 + 2p^{2n-7}(p-1) -4p +2p^{2n-4}(1-p^{-1}-2p^{-4}):=\mu_5(p, n),\\
\qquad\qquad\qquad\qquad \text{ if } \Gamma(G)=K_{p^{n-5}(p-1)}\cup (p+1)K_{p^{n-4}(p-1)}
\end{cases}
 \end{align*}
noting that $|V( K_{p^{n-4}(p^2-1)} \cup kpK_{p^{n-5}(p-1)}\cup (p-k)K_{p^{n-4}(p-1)})| =-p^{n-4}-p^{n-3}+2 p^{n-2}$, $|V( (kp+1)K_{p^{n-5}(p-1)}\cup (p+1-k)K_{p^{n-4}(p-1)})|=p^{n-2}-p^{n-5}$, $|V( K_{p^{n-4}(p^2-1)}\cup pK_{p^{n-4}(p-1)})|=-p^{n-4}-p^{n-3}+2 p^{n-2}$, $|V((p^2+p+1)K_{p^{n-5}(p-1)})|=p^{n-2}-p^{n-5}$ and $|V(K_{p^{n-5}(p-1)}\cup (p+1)K_{p^{n-4}(p-1)})|=p^{n-2}-p^{n-5}$.

\vspace{.3cm}

\noindent \textbf{Case 1.} $\Ecn(K_{|V(\Gamma (G))|})-\Ecn(\Gamma (G)) = \mu_1(p, n, k)$

In this case $p\geq 2$ and $n\geq 5$. For $p\geq 3$ and $n\geq 5$ we have $1\geq 3p^{-1}$, $p^{2n-8}>2$ and $p^{2n-7}\geq 2k$. Therefore, $\mu_1(p, n, k)>0$. 

If $p = 2$ and $n \geq 5$ then $	\mu_1(2, n, k)= (8kn-5) + (5k\times 2^{2n - 8} - 4k) + 7\times 2^{2n - 6} > 0$, since $k = 1,2$ and $2^{2n - 8} \geq 4$.

\vspace{.3cm}

\noindent \textbf{Case 2.} $\Ecn(K_{|V(\Gamma (G))|})-\Ecn(\Gamma (G)) = \mu_2(p, n, k)$ 

In this case $p\geq 2$ and $n\geq 5$. For $p\geq 5$ we have 

$$2p^{2n-4}(1-p^{-1})-4p^{2n-8}=2p^{2n-5}(p-1)-4p^{2n-8}\geq 0,$$ 
$$2p^{2n-6}(1-p^{-1})-4p=2p^{2n-7}(p-1)-4p > 0$$ 
and 
$$2kp^{2n-6}(1-3p^{-1})-4kp=2kp^{2n-7}(p-3)-4kp\geq 0.$$ 
Therefore, $\mu_2(p, n, k)>0$.
For $p=3$ and $n\geq 5$ we have 
$$\mu_2(3, n, k)=-16-8k+3^{2n-9}(352+16k) = (352\times 3^{2n-9} - 16) + (16k 3^{2n-9}-8k).$$
Since  $n\geq 5$, $3^{2n-9}\geq 3$ and so $\mu_2(3, n, k)>0$. 
For $p=2$ and $n\geq 5$ we have 
$$\mu_2(2, n, k)=-12-4k+2^{2n-8}k+18\times 2^{2n-8} = (18\times 2^{2n-8} -12) + (2^{2n-8}k-4k).$$
Since $n\geq 5$, $2^{2n-8}\geq 4$ and so $\mu_2(2, n, k)>0$.

\vspace{.3cm}

\noindent \textbf{Case 3.} $\Ecn(K_{|V(\Gamma (G))|})-\Ecn(\Gamma (G)) = \mu_3(p, n)$ 

In this case $p\geq 2$ and $n\geq 4$. Therefore, $3-5p^{-1} > 0$ and $p^{2n-6}+p^{2n-5}-2 > 0$. Hence,  $ \mu_3(p, n)>0$.

\vspace{.3cm}

\noindent \textbf{Case 4.} $\Ecn(K_{|V(\Gamma (G))|})-\Ecn(\Gamma (G)) = \mu_4(p, n)$ 


As shown in Case 2 of part (a), we have  $\mu_4(p, n)>0$ since $\mu_4(p, n) = \beta_2(p, n)$.

\vspace{.3cm}

\noindent \textbf{Case 5.} $\Ecn(K_{|V(\Gamma (G))|})-\Ecn(\Gamma (G)) = \mu_5(p, n)$

In this case $p\geq 2$ and $n\geq 5$. We have $2p^{2n-4}(1-p^{-1}-2p^{-4}) > 0$.  Since $p^{2n-9} > 1$ we have $4p^{2n-9} > 4$. Also,  $2p^{2n-7}(p-1) = (2p)p^{2n - 8}(p-1) > 4p$ since $2p \geq 4$ and $p^{2n - 8}(p-1) > p$. 
Therefore, $\mu_5(p, n)=>0.$

Thus, in all the cases $\Ecn(K_{|V(\Gamma (G))|})-\Ecn(\Gamma (G)) > 0$.
Hence, $\Gamma(G)$ is not CN-hyperenergetic.
\end{proof}

\begin{thm}
Let $G$ be a non-abelian  $p$-group of order $p^4$. Then $\Gamma(G)$ is not CN-hyperenergetic.
\end{thm}
\begin{proof}
By Theorem \ref{CNS2} and Corollary \ref{Thm 3.11}, we get
\begin{align*}
	\Ecn&(K_{|V(\Gamma (G))|})-\Ecn(\Gamma (G))\\
	&=\begin{cases} 
		2p(p^2(p-1)(p^2-1)-2):=\mu_1(p), &\text{ if } \Gamma(G)=(p+1)K_{p(p-1)}\\
		2p((p-1)+p^2(3p-5)):=\mu_2(p),  &\text{ if } \Gamma(G)=K_{p^2-1}\cup pK_{p-1}
\end{cases}
\end{align*}	
noting that $|V((p+1)K_{p(p-1)})|=p^3-p$ and $|V(K_{p^2-1}\cup pK_{p-1})|= 2p^2-p-1$.

\vspace{.3cm}

\noindent \textbf{Case 1.} $\Ecn(K_{|V(\Gamma (G))|})-\Ecn(\Gamma (G)) = \mu_1(p)$ 
	
Since $p\geq 2$ we have $p^2\geq 4$, $p-1\geq 1$, $p^2-1\geq 3$ and so $p^2(p-1)(p^2-1)\geq 12$. Therefore,  $p^2(p-1)(p^2-1)-2\geq 10$. Hence, Therefore, $\mu_1(p) >0$.

\vspace{.3cm}

\noindent \textbf{Case 2.} $\Ecn(K_{|V(\Gamma (G))|})-\Ecn(\Gamma (G)) = \mu_2(p)$

Since $p\geq 2$ we have $p-1>0$ and $3p-5>0$ and so $\mu_2(p) > 0$.
\end{proof}

\noindent \textbf{Concluding Remark:}  It is observed that the commuting conjugacy class graphs of all the groups considered in our paper are not CN-hyperenergetic. It may be interesting to find a finite group $G$ such that  $\Gamma(G)$ is CN-hyperenergetic or to prove that there is no finite group $G$ whose 
 $\Gamma(G)$ is CN-hyperenergetic.
 
\vspace{1cm}

\noindent \textbf{Acknowledgement:} 
The first author would like to thank DST for the INSPIRE Fellowship (IF200226).


\begin{thebibliography}{33}

\bibitem{ASG}
A. Alwardi,   N. D. Soner and  I. Gutman,  On the common-neighborhood energy of a graph, {\em Bulletin (Acad$\acute{\rm e}$mie Serbe Des Sciences Et Des Arts. Classe Des Sciences Math$\acute{\rm e}$matiques Et Naturelles.
Sciences Math$\acute{\rm e}$matiques}, {\bf 36}, 49--59, 2011.

\bibitem{BN-2021}
P. Bhowal and R. K. Nath, Spectral aspects of commuting conjugacy class graph of finite groups, {\em Algebraic Structures and Their Applications}, {\bf 8}(2), 67-118, 2021.

\bibitem{BN-2022}
P. Bhowal and R. K. Nath, Genus of commuting conjugacy class graph of certain finite groups, {\em Algebraic Structures and Their Applications}, {\bf 9}(1), 93-108, 2022.




\bibitem{FSN-2021}
W. N. T. Fasfous, R. Sharafdini and R. K. Nath, Common neighbourhood spectrum graphs of finite groups, {\em Algebra and Discrete Mathematics}, {\bf 32}(1), 33-48, 2021.

\bibitem{HLM}
M. Herzog, M. Longobardi and M. Maj, On a commuting graph on conjugacy classes of groups. {\em Communications in  Algebra}, {\bf 37}(10), 3369-3387, 2009.












\bibitem{MEFW-2016}
A. Mohammadian, A. Erfanian, D. G. M. Farrokhi and B. Wilkens, Triangle-free commuting conjugacy class graphs, {\em Journal of Group Theory}, {\bf 19}, 1049-1061, 2016.


\bibitem{NFDS-2021}
R. K. Nath, W. N. T. Fasfous, K. C. Das and Y. Shang, Common neighbourhood energy of commuting graphs of finite groups, {\em Symmetry}, {\bf 13}(9), 1651 (12 pages), 2021.


\bibitem{Salah-2020}
M. A. Salahshour, Commuting conjugacy class graph of $G$ when $\frac{G}{Z(G)}\cong D_{2n}$, {\em Mathematics Interdisciplinary Research}, {\bf 1}, 379-385, 2020.

\bibitem{SA-2020}
M. A. Salahshour and A. R. Ashrafi, Commuting conjugacy class graphs of finite groups, {\em Algebraic Structures and Their Applications}, {\bf 7}(2), 135-145, 2020.








\bibitem{SA-CA-2020}
M. A. Salahshour and A. R. Ashrafi, Commuting conjugacy class graph of finite CA-groups, {\em Khayyam Journal of Mathematics}, {\bf 6}(1), 108-118, 2020.



\bibitem{Walikar}
H. B. Walikar,  H. S. Ramane and  P. R. Hampiholi, On the energy of a graph,  {\em Graph connections}, Eds. R. Balakrishnan, H. M. Mulder and A. Vijayakumar,  Allied  publishers, New Delhi, 1999, pages 120--123.





\end{thebibliography}
\end{document}